\documentclass[
aps,%
12pt,%
final,%
notitlepage,%
oneside,%
onecolumn,
nobibnotes,%
nofootinbib,%
superscriptaddress,%
noshowpacs,%
centertags]%
{revtex4}

\usepackage{amsthm}
\usepackage[all]{xy}

\theoremstyle{plain}
\newtheorem*{Thm}{Theorem}
\newtheorem{Lem}{Lemma}

\begin{document}
\selectlanguage{english}

\title{Topology of the Spaces of Functions \\ with Prescribed Singularities on Surfaces}

\author{\firstname{E.~A.}~\surname{Kudryavtseva}}
\email{eakudr@mech.math.msu.su}
\affiliation{%
Moscow State University}%


\begin{abstract}
Let $M$ be a smooth connected orientable closed surface and $f_0\in C^\infty(M)$ a function having only critical points of the $A_\mu$-types, $\mu\in\NN$. Let $\cF=\cF(f_0)$ be the set of functions $f\in C^\infty(M)$ having the same types of local singularities as those of $f_0$. We describe the homotopy type of the space $\cF$, endowed with the $C^\infty$-topology, and its decomposition into orbits of the action of the group of ``left-right changings of coordinates''.

\end{abstract}

\maketitle

Let $M$ be a smooth connected orientable closed surface and $f_0\in C^\infty(M)$ a function having only critical points of the $A_\mu$-types, $\mu\in\NN$. Let $\cF=\cF(f_0)$ be the set of functions $f\in C^\infty(M)$ having the same types of local singularities as those of $f_0$. Denote by ${\mathcal D}^0(M)$ the identity component in the group ${\mathcal D}(M)={\rm Diff}^+(M)$ of orientation-preserving self-diffeomorphisms of $M$. The group ${\mathcal D}(\mathbb R)\times{\mathcal D}(M)$ acts on $\cF$ by ``left-right changings of coordinates''. We describe the homotopy type of the space $\cF$, endowed with the $C^\infty$-topology, and its decomposition into ${\mathcal D}(\mathbb R)\times{\mathcal D}^0(M)$-orbits. This result was announced in \cite {K:conf:Zlat,K:conf}. Similar result for a Morse function $f_0$ and $\chi(M)<0$ was obtained in \cite {KP3mn,KP4msb,KP4vest}.

Let us give a short historical overview, mostly for the case of a Morse function $f_0$ (see the paper \cite {KP3mn} and references therein). A.\,T.\ Fomenko posed the question (1997) whether the space $\cF$ is arcwise connected; it was answered affirmatively by the author \cite {Kud99} for $M=S^2,\RR P^2$, by S.\,V.\ Matveev \cite {Kud99} and H.\ Zieschang in the general case. Open ${\mathcal D}(\mathbb R)\times{\mathcal D}(M)$-orbits in $\cF$ were counted by V.\,I.\ Arnold \cite {Arn} and E.\,V.\ Kulinich (1998). Homotopy type of any ${\mathcal D}^0(M)$-orbit in $\cF$ was studied by S.\,I.\ Maksymenko \cite{Max} (when $f_0$ was allowed to have certain degenerate types singularities) and by the author \cite {KP3mn,KP4msb,KP4vest}. V.\,A.\ Vassiliev \cite {Vas} proved the parametric $h$-principle and studied cohomology of spaces of smooth ${\mathbb R}^N$-valued functions not having too complicated singularities on any smooth manifold $M$. However the 1-parameter $h$-principle fails for the spaces of Morse functions on some $M$ with $\dim M>5$ \cite {CL}.

\section{Main result} \label{sec:2}

For any function $f\in C^\infty(M)$, denote by $C_f$ the set of its ctitical points, and by $C_f^\triv$ the set of critical points of the $A_{2m}$-types, $m\in\NN$. Recall that, in a neighbourhood of such a point $x\in C_f$, there exist local coordinates $u,v$ such that $f=\eta(u^{2m+1}+v^2)+f(x)$ for some sign $\eta\in\{+,-\}$. The integer $\eta m$ will be called the {\em level} of the point $x$.

Denote by $C_f^\min$ and $C_f^\max$ (respectively $C_f^\saddle$) the set of critical points of $f$ of $A_{2m-1}$-types, $m\in\NN$, which are (respectively are not) points of local minima or local maxima. In a neighbourhood of such a point $x$, there exist local coordinates $u,v$ such that $f=\eta(u^{2m}\pm v^2)+f(x)$ where $\eta\in\{+,-\}$. The integer $\eta(m-1)$ will be called the {\em level} of the point $x$. The subset of degenerate critical points (i.e.\ those of non-zero levels) in $C_f^\extr:=C_f^\min\cup C_f^\max$ will be denoted by $\hat C_f^\extr$.

Suppose that an action of a group $G$ on a topological space $X$, a stratified \cite {Whi} orbifold $Y$ and a continuous surjection $\varkappa:X\to Y$ are given. If every $G$-orbit in $X$ is the full pre-image of a stratum from $Y$, we will say that $\varkappa$ {\it classifies} $G$-orbits, while $Y$ and $\varkappa$ are the {\it classifying space} and {\it map}.

The group $\cD(\RR)\times\cD(M)$ acts on $M\times\cF$ by the homeomorphisms $(x,f)\mapsto(h^{-1}(x),h_1^{-1}\circ f\circ h)$, $(h_1,h)\in\cD(\RR)\times\cD(M)$. 
Define the {\em evaluation functional} $\Eval:M\times\cF\to\RR$, $(x,f)\mapsto f(x)$, and
\begin{equation} \label {eq:1}
s:=\max\{0,\chi(M)+1\}>\chi(M).
\end{equation}

\begin{Thm} \label {thm:0}
For every function $f_0\in C^\infty(M)$ whose all critical points are of the $A_\mu$-types, $\mu\in\NN$, there exist smooth manifolds $\cB$ and $\cE$ and surjective submer\-sions 
$k:\cF\to\cB$, 
$\varkappa:M\times\cF\to\cE$, 
$\pi:\cE\to\cB$,
$\eps:\cE\to\RR$ such that
the diagram
\[
\xymatrix{
M\times\cF \ar[r]_{\varkappa} \ar@/^1.2pc/[rr]^\Eval \ar[d]_\Pr & \cE \ar[d]_{\pi} \ar[r]_{\eps} & \RR \\
\cF \ar[r]^{k} & \cB 
}\]
commutes, where $\Pr:M\times\cF\to\cF$ is the projection and $\dim\cB=2s+2|C_{f_0}^\triv|+|C_{f_0}^\extr|+|\hat C_{f_0}^\extr|+3|C_{f_0}^\saddle|=\dim\cE-2$. Moreover:

{\rm(a)} the maps $k,\varkappa$ are homotopy equivalences and classify $\cD^0(M)$- and $\cD(\RR)\times\cD^0(M)$-orbits in $\cF,M\times\cF$ for some strati\-fi\-ca\-tions on $\cB,\cE$ whose all strata are submanifolds; the map $\pi$ is a fibre bundle with fibres diffeomorphic to $M$;

{\rm(b)} the map $k$ (resp.\ $\varkappa$) induces a homotopy equivalence between every $\cD^0(M)$-invariant subset $B\subseteq\cF$ (resp.\ $E\subseteq M\times\cF$) and its image, e.g.\ between every orbit from item (a) and the corresponding stratum;

{\rm(c)} the group $\MCG(M)=\cD(M)/\cD^0(M)$ discre\-tely acts on $\cB,\cE$ by diffeo\-mor\-phisms preserving the strati\-fi\-ca\-tions from item (a) and the function $\eps$; the maps $p\circ k:\cF\to\cB':=\cB/\MCG(M)$ and $P\circ\varkappa:M\times\cF\to\cE':=\cE/\MCG(M)$ classify $\cD(M)$- and $\cD(\RR)\times\cD(M)$-orbits in $\cF$ and $M\times\cF$ for the induced stratifications on $\cB'$ and $\cE'$, where $p:\cB\to\cB'$ and $P:\cE\to\cE'$ are the projections.
\end{Thm}

Let us explain the term ``submersion'' in the case of functional spaces.
If $Q,R$ are smooth manifolds and $\cQ:=Q\times\cF$, denote by $C^\infty(R,\cQ)$ the preimage of $C^\infty(R,Q)\times C^\infty(R\times M)$ under the inclusion $C(R,\cQ)\hookrightarrow C(R,Q)\times C(R\times M)$, and by $C^\infty(\cQ,R)$ the set of maps inducing maps $C^\infty(\RR^n,\cQ)\to C^\infty(\RR^n,R)$ for all $\in\NN$. A map $p\in C^\infty(\cQ,R)$ will be called a {\em submersion} if, for any $q\in\cQ$, there exist a neighbourhood $U$ of the point $p(q)$ in $R$ and a map $\sigma\in C^\infty(U,\cQ)$ such that $p\circ\sigma=\id_U$.

\section{Constructing the classifying manifolds and maps} \label {sec:3}

Similarly to \cite {KP}, by a {\it framed function} on an oriented surface $M$ we will mean a pair $(f,\alpha)$ where $f\in C^\infty(M)$ has only the $A_\mu$-types local singularities and $\alpha$ is a closed 1--form on $M\setminus C_f^\extr$ such that (i) the 2-form $df\wedge\alpha$ has no zeros on ${M\setminus C_f}$ and defines a positive orientation, (ii) in a neighbourhood of every critical point $x\in C_f$ there exist local coordinates $u,v$ such that either $f=\eta(u^{2m+1}+v^2)+f(x)$ and $\alpha=\eta d(v-uv)$, or $f=\eta(u^{2m}-v^2)+f(x)$ and $\alpha=\eta d(uv)$, or $f=\eta(u^{2m}+v^2)+f(x)$ and $\alpha=\eta\varkappa_{f,x}\frac{udv-vdu}{u^2+v^2}$ where $\varkappa_{f,x}=\const>0$, $m\in\NN$, $\eta\in\{+,-\}$.

Denote by $\FF=\FF(f_0)$ the space of framed functions $(f,\alpha)$ such that $f\in\cF$. Endow this space with the $C^\infty$-topology \cite {KP}. Consider the right actions of $\cD(\RR)\times\cD(M)$ on $\FF$ and $M\times\FF$ by the homeomorphisms $(f,\alpha)\mapsto(h_1^{-1}\circ f\circ h,h^*\alpha)$ and $(x,f,\alpha)\mapsto(h^{-1}(x),h_1^{-1}\circ f\circ h,h^*\alpha)$, $(h_1,h)\in\cD(\RR)\times\cD(M)$.

Let $x_1,x_2,\dots\in M$ be pairwise distinct points. Denote by $\mathcal D^0_{r}(M)$ the identity component of the group $\cD_r(M):=\{h\in\cD(M)\mid h(x_i)=x_i,\ 1\le i\le r\}$, $r\in\ZZ_+$, whence $\cD_0(M)=\cD(M)$.

Define the classifying manifolds $\cB$ and $\cE$ as $\cB:=\cB_s$, $\cE:=\cE_s$, where $\cB_r$ and $\cE_r$ are the {\it universal moduli spaces} 
\[
\cB_r:= \FF/\cD^0_r(M), \quad 
\cE_r:= (M\times\FF)/\cD^0_r(M)
\]
{\it of framed functions} (resp. {\it framed functions with one marked point}) in $\cF$, $r\in\ZZ_+$. One shows similarly to \cite {KP3mn,KP4msb} that $\cB_r$ and $\cE_r$ are orbifolds of dimensions $\dim\cB_r=2r+2|C_{f_0}^\triv|+|C_{f_0}^\extr|+|\hat C_{f_0}^\extr|+3|C_{f_0}^\saddle|=\dim\cE_r-2$. For every group $\cG\in\{\cD^0(M),\cD(\RR)\times\cD^0(M)\}$, we endow $\cB_r$ and $\cE_r$ with the stratifications whose every stratum is the full preimage of a point under the projection $\cB_r\to\cF/\cG$ and $\cE_r\to(M\times\cF)/\cG$.

Due to the $\cD(M)$-equivariance of the projection $M\times\FF\to\FF$ and the $\cD(M)$-invariance of the evaluation functional
$M\times\FF\to\RR$, $(x,f,\alpha)\mapsto f(x)$, they induce some maps $\pi_r:\cE_r\to\cB_r$ and $\eps_r:\cE_r\to\RR$. Put $\pi=\pi_s$, $\eps=\eps_s$.

Similarly to \cite[Theorem 2.5]{KP} and \cite [Statement 5.3]{KP3mn}, one proves the following lemmata which readily imply the theorem.

\begin{Lem} \label {lem:1}
The projection $\Forg:\FF\to\cF$, $(f,\alpha)\mapsto f$, is a homotopy equivalence and has a homotopy inverse map $i:\cF\to\FF$ and corresponding homotopies that respect the projections $q:\cF\to\cF/\cD^0(M)$ and $q\circ\Forg:\FF\to\cF/\cD^0(M)$.
\end{Lem}

\begin{Lem} \label {lem:2}
If $r\ge s$ then $\cB_r$ is a smooth manifold, while the projection $\Ev_r:\FF\to\cB_r$ is a homotopy equivalence and has a homotopy inverse map $i_r:\cB_r\to\FF$ and corresponding homotopies that respect $\Ev_r$ {\rm(}whence $\Ev_r\circ i_r=\id_{\cB_r}${\rm)}.
\end{Lem}

Put $k_r=\Ev_r\circ i :\cF\to\cB_r$. One defines similarly $\varkappa_r$. Define the classifying maps $k=k_s$, $\varkappa=\varkappa_s$.

\section {Reducing to the case of Morse functions} \label {sec:4}

If $f_0$ is a Morse function and $s=0$, then the space $\cB$ from \S \ref {sec:3} coincides with the smooth stratified manifold $\tilde\cM$ (the {\it universal moduli space of framed Morse functions}) studied in \cite {KP3mn,KP4msb,KP4vest}. It happens that every $\cB_r$ and $\cE_r$ can be described in terms of Morse functions.

Recall that a function $f\in C^\infty(M)$ is said to be {\em Morse} if all its critical points are nondegenerate (i.e.\ have the $A_1$-type, cf.\ \S \ref{sec:2}). Denote by $\Morse(f_0)$ the space of Morse functions on $M$ having exactly $|C_{f_0}^\min|$ and $|C_{f_0}^\max|$ points of local minima and maxima and $|C_{f_0}^\saddle|$ saddle points.

A Morse function $f\in\Morse(f_0)$ will be called {\it $f_0$-labeled} if every its critical point $x\in C_f$ is labeled by an integer and, in the case when this integer does not vanish and $x\in C_f^\extr$, also by a 1-dimensional subspace $\ell_x\subset T_{x}M$, moreover $|C_{f_0}^\triv|$ of non-critical points of $f$ are labeled by non-zero integers in such a way that the level (cf.\ \S \ref{sec:2}) of every critical point of $f_0$ coincides with the integer label of the corresponding labeled point of $f$, for some bijections $C_{f_0}^\min\approx C_f^\min$, $C_{f_0}^\max\approx C_f^\max$, $C_{f_0}^\saddle\approx C_f^\saddle$ and a bijection between $C_{f_0}^\triv$ and the set of labeled non-critical points of $f$.

Denote by $\MMorse^*(f_0)$ the space of framed (cf.\ \S \ref {sec:3}) $f_0$-labeled Morse functions. It is not difficult to construct homeomorphisms
\begin{equation} \label {eq:2}
\cB_r\approx \MMorse^*(f_0)/\cD^0_r(M), \quad
\cE_r\approx (M\times\MMorse^*(f_0))/\cD^0_r(M), \qquad r\in\ZZ_+.
\end{equation}

\section {Relation with meromorphic functions and the configuration spaces} \label {sec:5}

Suppose that $M$ is either a sphere $S^2$ or a torus $T^2$. If $M=S^2$, denote by $\AA(f_0)$ the space of rational functions $R$ on the Riemann sphere $\overline\CC$ such that all poles of the 1-form $\omega=R(z)dz$ are simple and have real rezidues, being positive at $|C_{f_0}^\min|$ poles and negative at $|C_{f_0}^\max|$ poles. If $M=T^2$, denote by $\AA(f_0)$ the space of pairs $(\lambda,R)$ where $\lambda\in\CC$, $\Im\lambda>0$, and $R$ is a meromorphic function on the torus $T^2_\lambda=\CC/(\ZZ+\lambda\ZZ)$, whose poles are all simple, all periods of the meromorphic 1-form $\omega=R(z)dz$ are purely imaginary, and the residues are positive at $|C_{f_0}^\min|$ poles and negative at $|C_{f_0}^\max|$ poles.

Let $\AA_0(f_0)$ be the space of functions $R\in\AA(f_0)$ or pairs $(\lambda,R)\in\AA(f_0)$ such that $\omega=R(z)dz$ has only simple zeros.

Due to \cite [Proposition 3.4] {GruKri}, the assignment to a 1-form $\omega$ its poles and residues at them gives a bijection $\varphi:\AA(f_0)\stackrel{\approx}{\to} C(f_0)$, where $C(f_0)$ is the ``labeled configuration space'' consisting of $|C_{f_0}^\extr|$-points subsets of $M$ equipped by $|C_{f_0}^\min|$ positive and $|C_{f_0}^\max|$ negative real marks with zero total sum. Thus $\AA_0(f_0)$ is homeomorphic to the open subset $\varphi(\AA_0(f_0))\subseteq C(f_0)$ consisting of the ``labeled configurations'' that correspond to 1-forms $\omega$ without multiple zeros.

It is not difficult to derive from (\ref {eq:2}) with $r=s$ (cf.\ (\ref {eq:1}) and \cite[Remark 2.6]{KP}) that our manifold $\cB$ is homeomorphic to the space $\AA_0^*(f_0)$ of functions $R\in\AA_0(f_0)$ or pairs $(\lambda,R)\in\AA_0(f_0)$, marked by $f_0$-labels (cf.\ \S \ref{sec:4}) at zeros and poles of the 1-form $\omega=R(z)dz$ and at some other $|C_{f_0}^\triv|$ points, as well as by a ``vertical'' label consisting of (i) a real label and (ii) either a positive real label in the case of $|C_{f_0}^\triv|=|C_{f_0}^\saddle|=0$, or $|C_{f_0}^\extr|$ integral curves of the field $\ker(\Re\omega)$ separating the poles from other labeled points. Thus, the manifold $\cB\approx\AA_0^*(f_0)$ can be obtained from the ``labeled configuration subspace'' $\varphi(\AA_0(f_0))\subseteq C(f_0)$ by assigning the $f_0$-labels and the (topologically inessential) ``vertical'' label.

\begin{acknowledgments}
The author wishes to express gratitude to S.\,Yu.~Nemirovski for indicating the paper \cite {GruKri}. This work was done under the support of RFBR (grant \No~15-01-06302-a) and the programme ``Leading Scientific Schools of RF'' (grant NSh-7962.2016.1).
\end{acknowledgments}


\end{document}